\documentclass[10pt, english]{article}

\usepackage{latexsym, graphicx}
\usepackage{amsmath}
\usepackage{amsfonts}
\usepackage{amssymb}
\usepackage{hyperref}
\usepackage{babel}
\usepackage[babel=true]{csquotes}
\usepackage{geometry}
\usepackage{color}

\usepackage{amssymb,amsmath,amsfonts,mathrsfs,amsxtra,amsthm}
\usepackage{graphicx,color,subfig}
\usepackage{accents,xspace}
\usepackage{dsfont}
\usepackage{amscd}
\newcommand{\mypapersize}{
\setlength{\textwidth }{16cm}
\setlength{\textheight}{23.5cm}
\setlength{\oddsidemargin}{-0.14cm}
\setlength{\topmargin}{-1.6cm}
}

\mypapersize

\newtheorem{thm}{Theorem}[section]
\newtheorem*{thm*}{Theorem}

\newtheorem{oss}[thm]{Remark}

\usepackage{fancyhdr}
\begin{document}
\title{Part I: Rebuttal\, to \vspace{0.2cm}\\
\enquote{Uniform stabilization for the Timoshenko beam by a locally distributed
damping} \vspace{0.2cm}
}
\author{Fatiha Alabau-Boussouira\footnote{Laboratoire Jacques-Louis Lions Sorbonne Universit\'{e}, Universit\'{e} de Lorraine} 
}
\bibliographystyle{plain}
\date{July 26, 2023}
\maketitle {}
\tableofcontents
\begin{abstract}
A paper, entitled "Uniform stabilization for the Timoshenko beam by a locally distributed damping" was published in 2003, in the journal {\it Electronic Journal of Differential Equations}. Its title concerns exclusively its Section 3, devoted to the case of equal speeds of propagation and to its main theorem, namely Theorem 3.1. It states that the solutions of the Timoshenko system (see (1.3) in \cite{SW2003}) decays exponentially when the damping coefficient $b$ is locally distributed. The proof of Theorem 3.1 is crucially based on Lemma 3.6, which states the existence of a strict Lyapunov function along which the solutions of (1.3) decay when the speeds of propagation are equal.

This rebuttal shows the major gap and flaws in the proof of Lemma 3.6, which invalidate the proofs of Lemma 3.6 and Theorem 3.1. Lemma 3.6 is stated at the top of page 12.  The main part of its proof is given in the pages 12 and 13. In the last eight lines of page 13, eight inequalities are requested to hold together for the proof of Lemma 3.6. They don't appear in the statements of Lemma 3.6. The subsequent flaws come from the evidence that several of them are contradictory either between them or with claims in the title of the article. For instance, the seventh and eight inequalities stated in the last two lines of page 13 are basically contradictory. The first inequality among the eight ones implies that the damping coefficient $b$ is necessarily strictly positive on the domain, and thus is globally distributed. This contradicts the title of \cite{SW2003}.
Other examples of contradiction are described in this rebuttal. Therefore Lemma 3.6 remains unproved. As a consequence, the proof of Theorem 3.1 is false and the result claimed in the title of \cite{SW2003} is unproved.

We also point in this rebuttal other flaws, or gaps in the proofs of Theorem 2.2 related to strong stability and non uniform stability for the case of distinct speeds of propagation.

In \cite{AB2023II}, we correct and complete the proof of strong stability. We also correct, set up the missing functional frames, fill the gaps in the proof of non uniform stability in the cases of different speeds of propagation, and complete a missing argument in the proof of Theorem A in \cite{NRL1986} (see Remark \ref{ThmA}), the result of Theorem A being used in the paper \cite{SW2003} on which this rebuttal is mainly devoted.
\end{abstract}
\section{Introduction}
In this first part, we analyze the main result presented in the title of the paper ”Uniform stabilization for the Timoshenko beam by a locally distributed damping” published in 2003 in the journal Electronic Journal of Differential Equations written by Abdelaziz Soufyane and Ali Wehbe. This article will be refereed as \cite{SW2003} in the sequel.

The paper presents the Timoshenko system (1.1) introduced by Timoshenko in the reference appearing as [19] in \cite{SW2003}.  This model couples two wave equations through first order and zero order in space coupling terms. The coefficients $ \rho, K, I_{\rho}, E, I$ of the Timoshenko system are assumed to be constant and strictly positive. 
The purpose of \cite{SW2003} is to establish the well-posedness as well as stability properties of this system.
\vskip 2mm
We present here major flaws in the proof of the uniform stability of (1.3) in the case of equal speeds of propagation. We also present flaws and gaps in Theorem 2.2 concerning strong stability and non uniform stability when the speeds of propagation are distinct.

\vskip 2mm

We choose to present the rebuttal in a pedagogic way for an unaware public of readers, so that they can easily understand the major flaws with the concrete gaps, and the mathematical arguments explaining the flaws. 

\vskip 4mm
 
\section{Major gap and flaws in the proof of Lemma 3.6, Theorem 3.1 and other missing or contradictory arguments}

In \cite{SW2003}, the Timoshenko system is presented in (1.3). Its solutions are denoted by $w, \varphi$. When the full system is formulated as a first order system in time in \cite{SW2003}, the unknown  is denoted as $Y(t)=(w,\partial_tw, \varphi, \partial_t \varphi)^T$ , where the upperscript $T$ stands here for the transpose of a line vector to a column vector. 

 At page 13, the main conclusive steps of the proof of Lemma 3.6 are stated so that the authors can claim that if the speeds of propagation in (1.3) are equal, then the function $t \to \phi(Y(t))$ defined in (3.3) is a strict Lyapunov function, i.e. it satisfies the inequalities (3.1) and (3.2), where (3.2) is given at page 6 in \cite{SW2003} by:
 \enquote{
\begin{equation*}\label{R0}
\hspace{6cm} \partial_t\phi(Y(t)) \leqslant -d_2 \big|\big| Y(t) \big|\big|^2, \hspace{5.3cm}   (3.2)
\end{equation*}
}

where $d_2$ is a positive constant.  The assumption of equal speeds of propagation in (1.3) is given in following statement of Theorem 3.1 also at page 6 in \cite{SW2003}:

 \enquote{
\begin{equation*}
\dfrac{K}{\rho}=\dfrac{EI}{I_{\rho}}
\end{equation*}
}

The functions $t \to \phi(Y(t))$ defined in (3.3) and the function $t \to \phi_1(Y(t))$ (defined at page 8) of \cite{SW2003}, depend on the choices of four functions respectively denoted by $h,c, k, d$ which are supposed to admit derivatives of first oder for $h, k, d$ and up to second order for $c$ on the interval $[0,l]$ (the precise degree of smoothness is not explicitly stated in the paper). The functions $t \to \phi(Y(t))$ and the function $t \to \phi_1(Y(t))$ also depend on the choice of three positive and constant parameters $\varepsilon_i$ for $i=1,2, 3$. We also recall that $b$ is the damping function. 
 
 First the authors write an inequality concerning the function $t \to \phi(Y(t))$. They introduce in the estimates, positive constants $a_j$ for $j=1, \ldots, 7$ that can be chosen. We do not recopy exactly this inequality here but summarize it under the following form, to focus on the important points:
 \begin{multline*}\label{R1}
 \Big| \partial _{t}\phi (Y(t))\Big| \leqslant \int_0^l C^1_{b, h,c}(x)(\varphi_t(t,x))^2dx + \int_0^l C^2_{h,c,k}(x)(\varphi_x(t,x))^2dx +
 \int_0^l C^3_{c}(x)(\varphi(t,x))^2dx + \\ \int_0^l C^4 \big((\varphi - w_x)(t,x)\big)^2 dx +
  \int_0^l C^5_{k, d}(x)(w_x(t,x))^2dx +
 \int_0^l C^6_{k,d}(x)(w_t(t,x))^2dx + \\ C^7_{h(l)}(\varphi_x(l))^2 + C^8_{h(0)}(\varphi_x(0))^2 
 + C^9_{k(l)}(w_x(l))^2 + C^{10}_{k(0)}(w_x(0))^2.
 \end{multline*}
 Here, we choose to denote the coefficient functions $C^i_{\mbox{\tiny{subscript}}}$ where $i=1, \ldots 6$, which are space dependent only, by the appropriate subset of the functions $(h,k,c,d)$ on which they depend, and also for $C^1$ by the appropriate subset of the functions $(h,k,c,d)$ and $b$  (the damping coefficient). We denote by $C^i$
for $i=7, \ldots, 10$ the constants in front of the boundary values of $w_x$ and $\varphi_x$. Note that the constants also depend on the constant positive parameters to be adjusted: $\varepsilon_j$ for $j=1, 2, 3$ and $a_j$ for $j=1, \ldots, 7$, even though we do not include this dependence in the notation.

 The authors of \cite{SW2003} claim that the functions $h, c, k, d$ can be chosen, together with the parameters $\varepsilon_i$ for $i=1,2, 3$ (more specifically "sufficiently small"), and the constants  $a_j$ for $j=1, \ldots, 7$ in a judicious way to guarantee that $t \to \phi(Y(t))$ is a Lyapunov function, i.e. is such that the coefficient functions $C^i_{\mbox{\tiny{subscript}}}$ where $i=1, \ldots 6$ and the coefficients $C^i$
for $i=7, \ldots, 10$ respectively in front of $\varphi_t^2$, $\varphi_x^2$, $\varphi^2$, $\big( \varphi - w_x\big)^2$, $w_x^2$, $w_t^2$, and the boundary terms $(\varphi_x(l))^2$, $(\varphi_x(0))^2 $, $(w_x(l))^2$, $(w_x(0))^2$ are strictly negative on $(0,l)$. 

\vskip 2mm

For this property to hold, the authors state eight inequalities on the last eight lines at page 13. Let us present these eight inequalities and their purposes. The first coefficient function $C^1_{b, h,c}$ is splitted in $C^1_{b, h,c}=D^1_{b} + D^2_{h,c}$, and the two first inequalities request respectively that $D^1_{b} <0$ on $[0,l]$, and  $D^2_{h,c}<0$ on $[0,l]$. The third inequality requests $C^2_{h,c,k}<0$ on $[0,l]$. Note that a positive part of the coefficient $C^2_{h,c,k}$,  namely:
$$
c_{0}\Big(\dfrac{Kl^{2}\varepsilon _{1}\Big(2+\overline{c}^{2})}{2}\Big)\Big),
$$
has been forgotten by the authors, but this is a minor mistake. The fourth inequality requests that the constant coefficient $C^4<0$ (a minor typo is present but has no consequences). The fifth inequality requests that the coefficient function $C^6_{k,d}<0$ on $[0,l]$. The sixth inequality requests that the coefficient function $C^3_{c}<0$ on $[0,l]$. The seventh inequality requests that the coefficient function $C^5_{k, d}<0$ on $[0,l]$. The eighth inequality requests that all the constant coefficients $C^7_{h(l)}$, $C^8_{h(0)}$, $C^9_{k(l)}$ and $C^{10}_{k(0)}$ are strictly negative. 

\vskip 2mm

Let us precise now where the major gap and the flaw stand. The authors state that it is possible to choose the functions $h, c, k, d$, the positive parameters $\varepsilon_i$ for $i=1,2, 3$, and the positive constants $a_j$ for $j=1, \ldots, 7$ such that these eight inequalities are satisfied. This is asserted without any proof about their compatibilities. This is a major gap in the proof of Lemma 3.6. This gap also leads to the major flaws as we shall see. The authors also do not prove that these inequalities are compatible with their other assumptions, such as, for instance, their assumption that the support of the damping $b$ can be any non empty subset of the domain. This is claimed in the title of \cite{SW2003}. However, we shall also prove that one of these eight inequalities imply that the damping is necessarily globally distributed (a stronger assumption). Thus some of the eight inequalities are also not compatible with other assumptions of Lemma 3.6, and this of Theorem 3.1. 

Note that $\alpha$ and $\beta$ are positive constants defined in Lemma 3.4 at page 8 of \cite{SW2003} and that the constant $\alpha$ depends on $h$ in a nonlocal way (see \eqref{R15} below). 

We shall first prove that the seventh and eighth inequalities are incompatible, namely that they cannot be satisfied both.

The seventh inequality requests that $C^5_{k, d}<0$. Going back to \cite{SW2003}, one can check that the coefficient function $C^5_{k, d}$ (notation specific to this rebuttal) is given by:
$$
C^5_{k, d}=
$$
\enquote{
\begin{equation}\label{R2}
-\dfrac{K}{2}k_x+K\dfrac{a_5}{2}+\dfrac{K\overline{d}^{2}c_0}{2a_7}+
\dfrac{Kd^{2}}{2a_7}.
\end{equation}
}
\ \\
Thus the seventh inequality is equivalent to
\begin{equation}\label{R3}
K\dfrac{a_5}{2}+\dfrac{K\overline{d}^{2}c_0}{2a_7}+
\dfrac{Kd^{2}}{2a_7} < \dfrac{K}{2}k_x \mbox{ on } (0,l)
\end{equation}
Since the constants $K$, $a_5$, $c_0$ and $a_7$ are all strictly positive, we deduce that the seventh inequality requested in \cite{SW2003} implies that $k_x>0$ on $(0,l)$. Hence $k$ is a strictly increasing function on the interval $(0,l)$. Thus the seventh inequality requested in Soufyane-Wehbe's paper implies that
\begin{equation}\label{R4}
k(0)<k(l).
\end{equation}
Let us now turn to the eight inequality, which is given by
\begin{equation*}
\varepsilon _2<\min \Big(-\varepsilon _1 E I \dfrac{\rho }{K} h(l),\varepsilon
_1 E I\dfrac{\rho}{K} h(0),-\rho \varepsilon _3k(l),\rho \varepsilon_3k(0)\Big)
\end{equation*}
This implies in particular that
\begin{equation}\label{R5}
\varepsilon _2< -\rho \varepsilon _3k(l), \ \varepsilon _2< \rho \varepsilon_3k(0)
\end{equation}
Since $\varepsilon_i>0$ for $i=2,3$ and $\rho$ are all strictly positive, we deduce easily that the eight inequality requested in \cite{SW2003} implies
\begin{equation}\label{R6}
k(l)<0<k(0).
\end{equation}
Thus the two last inequalities at page 13, that is respectively the seventh and eighth inequalities imply \eqref{R4} and \eqref{R6} which are trivially incompatible. This is a major flaw, that shows that the authors did not prove that $t \to \phi(Y(t))$ is a Lyapunov function. As a consequence, the exponential stability property claimed by the authors in Theorem 3.1 is unproved. Other flaws are present in these eight inequalities. The first inequality among them is given by
$$
D^1_b=
$$
\enquote{\begin{equation*}
-b(x)+\varepsilon _1\dfrac{a_1b^2(x)}{2}+\varepsilon _1\dfrac{a_2b^2(x)}{
2}+\varepsilon _2\dfrac{I_{\rho }b^2(x)}{2a_6}<0
\end{equation*}}

Thus this inequality implies
\begin{equation}\label{R7}
\forall \ x \in \{y \in [0,l], b(y)=0\}, -b(x)+\varepsilon _1\dfrac{a_1b^2(x)}{2}+\varepsilon _1\dfrac{a_2b^2(x)}{
2}+\varepsilon _2\dfrac{I_{\rho }b^2(x)}{2a_6}=0 <0
\end{equation}
which does not hold true. Hence, the requested first inequality necessarily implies that
\begin{equation}\label{R8}
\{y \in [0,l], b(y)=0\}=\emptyset
\end{equation}
Since the authors assume that  $b$ is a continuous function on $[0,l]$ and such that $b \geqslant 0$ on $[0,l]$ (see page 2), this implies that
\begin{equation}\label{R9}
b>0 \mbox{ on } [0,l],
\end{equation}
so that 
\begin{equation}\label{R10}
\mbox{ There exists a constant }\ b_->0 \mbox{ such that } b \geqslant b_- \mbox{ on } [0,l].
\end{equation}
Thus the first requested inequality at page 13 in \cite{SW2003} implies that the damping is globally distributed and not locally distributed as stated in the title of the article at least for Theorem 3.1 and Lemma 3.6. 

Let us make the observation that the statement of Lemma 3.6 includes the following hypotheses on the functions
$c,h,k,d$ 
\enquote{
\begin{equation*}
\begin{cases}
-\dfrac{h_{x}(x)}{2}+c(x) <0\quad\forall \  x \in ]0,l[\backslash]b_0,b_1[
\\
\dfrac{h_{x}(x)}{2}+c(x) >0\quad \forall \  x \in ]0,l[ \\
\dfrac{k_{x}(x)}{2}+d(x) > 0\quad \forall \ x \in ]0,l[,
\end{cases}
\end{equation*}
}
\ \\
However, the second inequality of the eight requested inequalities at page 13, which is equivalent to the request
$$
D^2_{h,c} = \varepsilon_1 I_{\rho}\big(-\dfrac{h_{x}}{2}+c\big) + \dfrac{\varepsilon_1\alpha}{2a_3} + \dfrac{\varepsilon_1\beta c^2}{2a_4} + \varepsilon_2<0
$$
 should hold on all $(0,l)$, so that in particular
\begin{equation}\label{R11}
-\dfrac{h_{x}}{2}+c <0, \mbox{ on } [b_0,b_1]
\end{equation}
is also requested. Thus, adding this additional assumption in the statement of Lemma 3.6, one has
\begin{equation}\label{R12}
\begin{cases}
-\dfrac{h_{x}}{2}+c <0, \mbox{ on } (0,l),
\\
\dfrac{h_{x}}{2}+c>0, \mbox{ on } (0,l),
\end{cases}
\end{equation}
which imply
\begin{equation}\label{R13}
h_x >0, \mbox{ on } (0,l).
\end{equation}
But this is contradictory once again with the eighth requested inequality (last line) at page 13, since it implies 
\begin{equation}\label{R14}
h(l)<0<h(0).
\end{equation}
Thus several flaws are both in the statement and the proof of Lemma 3.6. 

\begin{oss}\label{Rk4}
Let us make more comments on this section.
One can note that the sign of absolute value present in the estimate of $\Big| \partial _{t}\phi (Y(t))\Big|$ should be withdrawn since the authors state inequalities in such a way that $\Big| \partial _{t}\phi (Y(t))\Big|$  is smaller or equal to a negative term.
\ \\
In Lemma 3.3, the assumptions on the functions $h$ and $c$ are not needed, whereas one has to assume that the initial data for $Y$ solution of (1.3) are in the domain of $L$ to justify rigorously the computations. Moreover, one can note that these assumptions are contradictory with the property $h_x>0$ on $(0,l)$.
The assumptions on the functions $h,c, k, d$ are more generally not used in Lemma 3.4 and 3.5. The requested assumptions arise from the eight inequalities stated in the last lines of page 13. The assumption of equal speeds of propagation is used in Lemma 3.5 but is not in its statement.
\ \\
In Lemma 3.4, the definition of $\alpha$ and $\beta$ should be replaced due to mathematical coherence, by
\begin{equation}\label{R15}
\alpha =| (h_{x}\partial _{x}(-\partial _{xx})^{-1}-hI_{L^2})|
_{\mathcal{L}(L^{2}(0,l))}\,,\quad \beta =| \partial _{x}(-\partial _{xx})^{-1}|
_{\mathcal{L}(L^{2}(0,l))}
\end{equation}
Note also that, $\varphi$ should be replaced by $w$ in the left hand side of the last inequality at page 10. This typo has no consequences on the sequel. Some other typos exist but they have no consequences.
\end{oss}

\section{Flaw and gap in the proof of Theorem 2.2 in Section 2. of \cite{SW2003}}

\subsection{A first flaw in the proof of the strong stability result in Theorem 2.2  of \cite{SW2003}}
The proof of the statement that $L$ has no pure imaginary eigenvalues in Theorem 2.2 is given at page 4.  The proof is by contradiction. The authors assume that $L$ has a pure imaginary eigenvalue $i\omega$, where $\omega \in \mathbb{R}^{\ast}$, with a corresponding non vanishing eigenvector $U_1$ and claim that it leads to a contradiction, namely that $U_1 \equiv 0$. The property that $U_1$ is an eigenvector of $L$ associated to the eigenvalue $i\omega$ with $U_1=(u, u_2,v,v_2)^T \in D(L)$ implies that $(u,v)$ is a non vanishing solution of the coupled second order system stated in \cite{SW2003}.

\enquote{
\begin{equation*}\label{R16}
\begin{cases}
Ku_{xx} -Kv_x=-\rho\, \omega^2 u, \quad \mbox{in } (0,l),\\
EIv_{xx} +Ku_x -Kv-ib\omega v=-I_{\rho} \omega^2 v, \quad \mbox{in } (0,l), \hspace{7.5cm}   (2.5)\\
u(0)=u(l)=v(0)=v(l)=0.
\end{cases}
\end{equation*}
}

\noindent One has furthermore the relations
\begin{equation*}\label{R16b}
u_2=i\omega u \quad \mbox{in } (0,l), v_2=i \omega v  \quad \mbox{in } (0,l).
\end{equation*}

Note also that both $u$ and $v$ have complex-valued functions, and that one cannot perform computations assuming that they are real-valued without proving this property. This is due to the presence of the damping term in the second equation of (2.5) (namely the term $ib\omega v$). The authors start their proof by multiplying the second equation of (2.5) by $v$, and then assume that $v_x^2=|v_x|^2$, $v^2=|v|^2$, $u_x v \in \mathbb{R}$ which do not hold for complex-valued functions without further arguments of proof. This leads to the a first main flaw in the proof. 
\begin{oss}\label{Rk0}
Let us explain more precisely the flaw, due to the fact that $u$ and $v$ are complex-valued functions, one has multiply the second equation of (2.5) by $\overline{v}$ instead of $v$, integrate by part over the interval $[0,l]$ and perform an integration by parts on the first term involving $\int_0^l EI\partial_{xx}v\overline{v}dx=-\int_0^l EI \Big|v_x\Big|^2dx$, due to the homogeneous Dirichlet boundary condition on $v$ at $x=0$ and $x=l$. This leads to the equality 

\begin{equation}\label{R17}
EI\int_0^l \Big|\partial_xv\Big|^2 dx + \int_0^l \Big(K-I_{\rho}\omega^2\Big)|v|^2dx +i\omega \int_0^l b|v|^2dx -K \int_0^l \partial_xu \overline{v}dx=0.
\end{equation}

The first two integral terms are real numbers.  If $\int_0^l b|v|^2dx \neq 0$, the third term is a purely imaginary number, whereas one cannot assert that the last term in \eqref{R17} has a vanishing imaginary part nor a vanishing real part. Hence, the statement written in \cite{SW2003}, next to (2.5) that:
\ \\
\enquote{
\begin{equation*}
\int_{0}^{l}(EI| \partial _{x}v| ^{2}-K\partial _{x}u.v+(K-I_{\rho }\omega
^{2})| v| ^{2})dx=0\;\text{\ and\ \ }\int_{0}^{l}\omega b(x)|v|^{2}dx=0,
\end{equation*}
}
\ \\
is unproved due to the above flaw. This also has consequences on the proof that $(u,v)$ satisfies the system (2.6) in \cite{SW2003}, namely
\enquote{
\begin{equation*}\label{R23}
\begin{cases}
Ku_{xx} -Kv_x=-\rho\, \omega^2 u, \quad \mbox{in } (0,l),\\
EIv_{xx} +Ku_x -Kv=-I_{\rho} \omega^2 v, \quad \mbox{in } (0,l), \hspace{8.5cm} (2.6)\\
v=u=0, \quad \mbox{in } (b_0,b_1),\\
u(0)=u(l)=v(0)=v(l)=0.
\end{cases}
\end{equation*}
}

\end{oss}
\begin{oss}\label{Rk1}
If the damping coefficient is assumed to be vanishing on $[0,l]$, then one can assume that $u$ and $v$ are real-valued. However, in this case the energy of the solutions of (1.3) is conserved through time so that the system is not strongly stable. 
\end{oss}

\subsection{A gap in the proof of strong stability in Theorem 2.2 of \cite{SW2003}}
In \cite{SW2003}, the authors conclude quickly without any proof that $u=v=0$ is the solution of \eqref{R23}. Effectively, $u=v=0$ is {\it a} solution of \eqref{R23}, the conclusive argument being to prove that it is the only solution. Whenever rigorously proved, the property is called a "unique continuation property". 
\section{On gaps in the proof of non-uniform stability in Theorem 2.2 in case of distinct speeds of propagation}

The second statement in Theorem 2.2 concerns the property of non-uniform stability when 
\begin{equation*}
\dfrac{K}{\rho} \neq \dfrac{EI}{I_{\rho}}
\end{equation*}
Under this hypothesis, Soufyane and Wehbe can make use of an approach presented in [6] (referenced in the present rebuttal as \cite{NRL1986}). They first set $L_1=L-D$ where $D$ is a suitable compact operator in the energy space. The difference between $L$ and $L_1$ being compact, the $\mathcal{C}^0$ semigroups $(e^{tL})_{t \geqslant 0}$ and $(e^{tL_1})_{t \geqslant 0}$ generated respectively by $L$ and $L_1$ have the same essential spectral radius, i.e. $r_e(e^{tL})=r_e(e^{tL_1})$ for all $t \geqslant 0$.  The coupled system corresponding to the new first order system $\widehat{Y}'(t)=L_1\widehat{Y}(t)$,
where $\widehat{Y}(t)=(u,u_t,v,v_t)^T$, is similar to the system (1.3) up to the suppression of the term $-K\varphi$ in the right hand side of the second equation in (1.3), $v$ replacing $\varphi$ in the new system. More precisely $(u,v)$ solves now the coupled system:
\begin{equation*}\label{R4.1}
\begin{cases}
u_{tt}=\dfrac{K}{\rho}(u_{x}-v)_{x},  \quad t\geqslant 0, x \in (0,l), \\[1em]
v_{tt}=\dfrac{EI}{I_{\rho }}v_{xx}+\dfrac{K}{I_{\rho }}u_{x}-\dfrac{b}{I_{\rho }}v_{t}, \quad t\geqslant 0, x \in (0,l), \\
u(0,t)=u(l,t)=0;\quad v(0,t)=v(l,t)=0, \quad t \geqslant 0.
\end{cases} 
\end{equation*}
The authors follows [6] (that is \cite{NRL1986} in this rebuttal) using the Riemann invariants $(p,\varphi, q, \psi)^T$ associated to the first order system $Y'(t)=L_1Y(t)$ 

\begin{oss}\label{Rk3}
In \cite{SW2003}, one should replace the constant diagonal matrix $K$ of rank $4$ by $\widehat{K}$ (for instance) in (2.7) and in the top of page 5, the notation $K$ being already used by the authors in (1.3).
In a similar way, one should replace the notation $\varphi$ for the second Riemann invariant by the notation $\hat{\varphi}$, the notation $\varphi$ being already used in (1.3). 
\end{oss}
The Riemann invariants are given at page 4 of \cite{SW2003} by:

\enquote{
\begin{equation*}\label{R4.2}
\begin{cases}
p=-\sqrt{\dfrac{K}{\rho}}u_x + u_t, \quad q= \sqrt{\dfrac{K}{\rho}}u_x + u_t,\\
\hat{\varphi}=-\sqrt{\dfrac{EI}{I_{\rho}}}v_x + v_t, \quad \psi= \sqrt{\dfrac{EI}{I_{\rho}}}v_x + v_t.
\end{cases} 
\end{equation*}
}

This leads to the first order system given in \cite{SW2003}:
\begin{equation*}
\partial_tZ + \widehat{K}\partial_xZ +C Z=0
\end{equation*}
subjected to the homogeneous Dirichlet boundary conditions for $(p+q)$ and $\hat{\varphi} + \psi$ given at the top of page 5 in \cite{SW2003}, where 
 $Z=(p, \hat{\varphi}, q, \psi)^T$, $\widehat{K}$ s a diagonal square matrix of order $4$ with constant coefficients depending  on the physical constants involved in (1.3), and where $C$ is a square matrix of order $4$, whose coefficients depend on the damping coefficient $b$ and on the physical constants involved in (1.3). 
 
 Note that a typo is present in \cite{SW2003}, the sign minus in front of the coefficient $b$ in all the elements of the matrix $C$ has to be changed to a sign plus.

Let us define the operator:
\begin{equation}\label{R4.3}
S_1:=-(\widehat{K}\partial_x +C)
\end{equation}
associated to homogeneous Dirichlet boundary conditions for $(p+q)$ and $\hat{\varphi} + \psi$ at both ends $x=0$ and $x=l$:
\begin{equation}\label{R4.4}
(p+q)(0)=(p+q)(l)=(\hat{\varphi} + \psi)(0)=(\hat{\varphi} + \psi)(l)=0.
\end{equation}
Let us denote by $C_0$ the square matrix of order $4$ corresponding to the diagonal of $C$.
 The authors of \cite{SW2003} then rely on a result by Neves, Ribeiro and Lopes in [6] (that is \cite{NRL1986} in this rebuttal), valid when the positive as well as the negative eigenvalues of the diagonal matrix $\widehat{K}$ are distinct, which holds true in \cite{SW2003} thanks to the hypothesis of distinct speeds of propagation, to claim that the essential spectral radius of of $L_1$ is equal to the essential spectral radius of the operator
 $L_2=-(\widehat{K}\partial_x + C_0)$ associated to the boundary conditions \eqref{R4.4}.
 \begin{oss}\label{Rk4b}
 To understand the gaps in \cite{SW2003}, let us precise the framework introduced by Neves-Ribeiro-Lopes in [6] (\cite{NRL1986} in this rebuttal). Let $M$ be a given diagonal matrix with eigenvalues ordered as strictly $N$ positive ones and $n-N$ strictly negative ones, where the matrices $D, F, G, \widehat{E}$ have the appropriate lines and columns and $n\,, N \in \mathbb{N}^{\ast}$, with $n >N$.
 Neves-Ribeiro-Lopes's article [6] (\cite{NRL1986} in this rebuttal) compares the essential spectral radius of first order hyperbolic systems involving $M$ and a zero order operator $C$, subjected to certain boundary conditions, either dynamic or static, under a different frame of boundary conditions than the one present in \cite{SW2003}: namely in the autonomous case, [6] (\cite{NRL1986} in this rebuttal) compares a first hyperbolic system in a block form written as 
 \begin{equation}\label{R4.5}
 \partial_t(u,v)^T + M\partial_x(u,v)^T + C(u,v)^T=0, \mbox{ on } [0,l],
 \end{equation}
 subjected to the dynamic boundary condition
 \begin{equation}\label{R4.6}
 \dfrac{d}{dt}\Big(v(t,l)-Du(t,l)
 \Big)=Fu(t,l) + Gv(t,l), 
 \end{equation}
and a boundary condition of the form
\begin{equation}\label{R4.7}
 u(t,0)=\widehat{E}v(t,0)\
 \end{equation}
\noindent to the first hyperbolic system in a block form written as
\begin{equation}\label{R4.8}
 \partial_t(u,v)^T + M\partial_x(u,v)^T + C_0(u,v)^T=0, \mbox{ on } [0,l],
 \end{equation}
\begin{equation}\label{R4.9}
 v(t,l)=Du(t,l), 
 \end{equation}
and \eqref{R4.7}, where $C_0$ the square matrix corresponding to the diagonal of $C$ in [6] (\cite{NRL1986} in this rebuttal). The column vector functions $u$ and $v$ are respectively of size $N$ and $n-N$.

\end{oss}

Under appropriate hypotheses the authors of [6] (\cite{NRL1986} in this rebuttal) prove that the essential spectral radius $r_e^C$ of the semigroup generated by the operator $-(M\partial_x +C)$ subjected \eqref{R4.6}-\eqref{R4.7} is equal to the essential spectral radius $r_e^{C_{0}}$ of the semigroup generated by $-(M\partial_x +C_0)$ subjected to \eqref{R4.7} and  \eqref{R4.9} (see problem (II.4) in \cite{NRL1986}). This is this property which is claimed to be used at page 5 of  in \cite{SW2003} as briefly mentioned just after the definition of the matrix $C$. Indeed we shall see that this is not sufficient to prove the claimed equality in \cite{SW2003}.

Note that the authors of \cite{NRL1986} prove a much more general result, also considering that the coefficients of the matrices $M$, and $C$ can depend on time and space and the coefficients of the matrices $\widehat{E}$, $D$, $F$, and $G$ can depend on time. 

\begin{oss}\label{ThmA}
Note that the static boundary conditions \eqref{R4.4} (which are given in the second line at page 5 in \cite{SW2003}) correspond to the boundary conditions \eqref{R4.7} and  \eqref{R4.9} with the appropriate choices of matrices $\widehat{E}$, $D$, $F$ and $G$. Note that it seems that an argument is missing in \cite{NRL1986} for the complete proof in Theorem A of the equality of the essential radius of the  semigroups generated by the operator $A$ and the operator $A_4$. This equality being used in \cite{SW2003}, we complete the proof of Theorem A of \cite{NRL1986} by a more general result stated in Proposition 4.4 (and a more general formulation in Proposition 4.6), which is proved in \cite{AB2023II}. This  result allows to complete the proof of Theorem A in a particular case of application. The result has also to be used a second time to fill a gap in \cite{SW2003}.
\end{oss}

\vskip 2mm

 Let us now describe the gaps in \cite{SW2003} in the proof of non uniform stability when the speeds of propagation in (1.3) are assumed to be distinct. Let us denote by $S_2$ the operator $-(\widehat{K}\partial_x +C)$ subjected to the boundary conditions \eqref{R4.6}-\eqref{R4.7}. Note now, that $S_1$ is the operator $-(\widehat{K}\partial_x +C)$ subjected to the boundary conditions \eqref{R4.7} and \eqref{R4.9}. Thanks to Neves-Ribeiro-Lopes [6] (\cite{NRL1986} in this rebuttal), $r_e(e^{tS_2})=r_e(e^{tL_2})$.
 It is not proved in \cite{SW2003} that $r_e(e^{tS_1})=r_e(e^{tS_2})$. This is a first gap in the proof. The second gap is that Soufyane and Wehbe do not prove in \cite{SW2003} that $r_e(e^{tL_1})=r_e(e^{tS_1})$. Note also the mathematical links between these different operators are missing. One has also to define them precisely. For this, one has also to define the domains of the operators $S_1$ and $S_2$. The domains of $S_1$ and $S_2$ are not defined in \cite{SW2003}. The definitions of these domains will also take into account the static boundary conditions, the dynamic ones having to be considered in a slightly different way. 
 
 \vskip 2mm
 In \cite{SW2003}, the authors use once again \cite{NRL1986} to write at page 5 that:
 
 \enquote{
 $$
 r_e(e^{tL_2})=\exp{(\alpha\, t)} \mbox{ where}
 $$
 }
 
  \enquote{
 $$
 \alpha=\sup\{\Re(\lambda), \lambda \in \sigma(L_2)\}
 $$
 }
 
 Then, they compute the punctual spectrum of $L_2$ and prove that 
 $$
 \alpha=\sup\{\Re(\lambda), \lambda \in \sigma_p(L_2)\}
 $$
 where $\sigma_p(L_2)$ denotes the punctual spectrum of $L_2$, that is the set of the eigenvalues of $L_2$. However in the general case, the full spectrum of an unbounded 
 operator, is strictly larger and not reduced to its punctual spectrum. Hence further mathematical arguments are missing in \cite{SW2003} to
 assert that:
 \begin{equation}\label{Punctual}
 \sup\{\Re(\lambda), \lambda \in \sigma_p(L_2)\}=\sup\{\Re(\lambda), \lambda \in \sigma(L_2)\}
 \end{equation}
 \vskip 2mm
 
In \cite{AB2023II}, we fill these gaps and correct the proof of non uniform stability given in \cite{SW2003}.
  
\section{Description of the presentation, open access, submission history}

The article is available in open access at the link: 
\vskip 0.2 cm
\url{https://ejde.math.txstate.edu/Volumes/2003/29/soufyane.pdf}
\vskip 0.2 cm
The pdf of the article starts at page 1 and ends at page 14.
\vskip 0.2 cm
The history of submission is indicated at the bottom of page 1 as follows:
\vskip 0.2cm 
\enquote{Submitted September 2, 2002. Published March 16, 2003}
\vskip 0.2 cm

Acknowledgments at page 223 mention that:
“The authors wish to thank the anonymous referee for his or her valuable suggestions. A. Soufyane would like to thank Dr. A. Benabdallah for her help.”

\vskip 0.2 cm

The publisher and main sponsor of the Electronic Journal of Differential Equations is the Department of Mathematics,Texas State University. 

\vskip 0.2 cm

The full journal is on open access and provides "rapid dissemination of high quality research in mathematics" as mentioned on the website of the journal.

\section{Data on dissemination through citations}
At the time this rebuttal was written the citations of \cite{SW2003} co-authored by Soufyane and Wehbe in some of the existing databases or platforms are 
\begin{center}
\begin{tabular}{|l|r|}
  \hline
  Google scholar citations & 203 \\
  \hline
  Web of Science citations & 126 \\
  \hline
  MathSciNet citations & 115 \\
  \hline
  zbMATH citations & 49 \\
  \hline
\end{tabular}
\end{center}

\begin{oss}\label{WOS}

The citations may vary. Concerning  Web of Science citations, it appears that the citations are splitted into 4 groups of citations: a group of 91 citations, a group of 33 citations, 1 group of 1 citation, 1 group of 1 citation. Last June, we could count 128 Web of Science citations, also splitted in 4 groups of citations with apparently non redundancy (94 citations, 32 citations, 1 citation, 1 citation).

\end{oss}

\subsection{Some description of 115 citations in the MathSciNet databasis}

\subsubsection{Distribution by Publishers and peer-reviewed journals}

\begin{itemize}

\item Publisher Springer : 31 citations distributed as follows

\begin{itemize}

\item 10 citations in Zeitschrift f\"ur Angewandte Mathematik und Physik

\item 4 citations in Acta Mathematica Scientia. Series B. English Edition

\item 3 citations in Applied Mathematics and Optimization

\item 3 citations in Journal of Dynamical and Control Systems

\item 3 citations in Mediterranean Journal of Mathematics

\item 1 citation in Journal of Evolution Equations

\item 1 citation in Acta Applicandae Mathematicae

\item 1 citation in Bulletin of the Malaysian Mathematical Sciences Society

\item 1 citation in Journal of Dynamics and Differential Equation

\item 1 citation in Archiv der Mathematik (Basel)

\item 1 citation in Boundary Value Problem

\item 1 citation in Nonlinear Differential Equations and Applications

\item 1 citation in Qualitative Theory of Dynamical Systems

\end{itemize}

\item Publisher Elsevier Inc.: 17 citations distributed as follows

\begin{itemize}

\item 11 citations in Journal of Mathematical Analysis and Applications

\item 2 citations in Nonlinear Analysis Theory Methods \& Applications

\item 2 citations in Nonlinear Analysis Real World Applications

\item 1 citation in Automatica Journal of the International Federation of Automatic Control

\item 1 citation in Journal of Differential Equations

\end{itemize}

\item Publisher Taylor \& Francis: 14 citations distributed as follows

\begin{itemize}

\item 14 citations in Applicable Analysis

\end{itemize}

\item Publisher Wiley: 15 citations distributed as follows

\begin{itemize}

\item 13 citations in Mathematical Methods in the Applied Sciences

\item 2 citations in Zeitschrift f\"ur Angewandte Mathematik und Mechanik

\end{itemize}

\item Publisher : American Institute of Mathematical Sciences 8 citations distributed as follows

\begin{itemize}

\item 2 citations in Evolution Equations and Control Theory

\item 2 citations in Communications on Pure and Applied Analysis

\item 1 citation in Discrete and Continuous Dynamical Systems Series S

\item 1 citation in Discrete and Continuous Dynamical Systems Series A

\item 1 citation in Discrete and Continuous Dynamical Systems Series B

\item 1 citation in Mathematical Control and Related Fields

\end{itemize}

\item Publisher :  American Institute of Physics 5 citations distributed as follows

\begin{itemize}

\item 5 citations in Journal of Mathematical Physics

\end{itemize}

\item Publisher : Dynamic Publishers Inc.  4 citations distributed as follows

\begin{itemize}

\item 4 citations in Dynamic Systems and Applications

\end{itemize}

\item Publisher : IOS Press  3 citations distributed as follows

\begin{itemize}

\item 3 citations in Asymptotic Analysis

\end{itemize}

\item Publisher : Oxford Academic Press  3 citations distributed as follows

\begin{itemize}

\item 2 citations in IMA Journal of Applied Mathematics

\item 1 citation in The Quarterly Journal of Mechanics and Applied Mathematics

\end{itemize}

\item Publisher : Khayyam  3 citations distributed as follows

\begin{itemize}

\item 2 citations in Differential and Integral Equations

\item 1 citation in Advances in Differential Equations

\end{itemize}

\item Publisher :  Acad\'emie des Sciences Paris 2 citations distributed as follows

\begin{itemize}

\item 2 citations in Comptes Rendus Math\'ematique

\end{itemize}

\item Publisher : Texas State Univ.--San Marcos, Dept. Math 2 citations distributed as follows

\begin{itemize}

\item 2 citations in Electronic Journal of Differential Equations

\end{itemize}

\item Publisher :  Sage 2 citations distributed as follows

\begin{itemize}

\item 2 citations in Mathematics and Mechanics of Solids

\end{itemize}

\item Publisher :  Osaka Univ. Osaka City Univ 1 citation distributed as follows

\begin{itemize}

\item 1 citation in Osaka Journal of Mathematics

\end{itemize}

\item Publisher :  Rocky Mountain Math. Consortium 1 citation distributed as follows

\begin{itemize}

\item 1 citation in Journal of Integral Equations and Applications

\end{itemize}

\item Publisher :  EDP Science 1 citation distributed as follows

\begin{itemize}

\item 1 citation in Mathematical Modelling of Natural Phenomena

\end{itemize}

\item Publisher :  De Gruyter 1 citation distributed as follows

\begin{itemize}

\item 1 citation in Advances in Nonlinear Analysis

\end{itemize}

\item Publisher :  Institute of Electrical and Electronics Engineers (IEEE) 1 citation distributed as follows

\begin{itemize}

\item 1 citation in IEEE Transactions on Automatic Control

\end{itemize}

\end{itemize}

\subsubsection{Citations of \cite{SW2003} by authors}

\begin{itemize}

\item 23 citations by Salim A. Messaoudi, collaboration distance with Abdelaziz Soufyane =1, collaboration distance with Ali Wehbe=1

\item 17 citations by Muhammad I. Mustafa, collaboration distance with Abdelaziz Soufyane=2, collaboration distance with Ali Wehbe=2

\item 10 citations by Abdelaziz Soufyane

\item 10 citations by Nasser-eddine Tatar, collaboration distance with Abdelaziz Soufyane=2, collaboration distance with Ali Wehbe=2

\item 10 citations by Dilberto da Silva Almeida J\'unior, collaboration distance with Abdelaziz Soufyane=1, collaboration distance with Ali Wehbe=2

\item 9 citations by Belkacem Said-Houari, collaboration distance with Abdelaziz Soufyane=1, collaboration distance with Ali Wehbe=2

\item 8 citations by Baowei Feng, collaboration distance with Abdelaziz Soufyane=1, collaboration distance with Ali Wehbe=2

\item 7 citations by Tijani Abdul-Aziz Apalara, collaboration distance with Abdelaziz Soufyane=1, collaboration distance with Ali Wehbe=2

\item 7 citations by Ali Wehbe

\item 6 citations by Aissa Guesmia, collaboration distance with Abdelaziz Soufyane =1, collaboration distance with Ali Wehbe=1

\item 6 citations by Ramos, Anderson de Jesus Ara\'ujo, collaboration distance with Abdelaziz Soufyane=1, collaboration distance with Ali Wehbe=2

\item 5 citations by Mauro de Lima Santos, collaboration distance with Abdelaziz Soufyane=1, collaboration distance with Ali Wehbe=2

\item 4 citations by Mohammad M. Kafini, collaboration distance with Abdelaziz Soufyane=2, collaboration distance with Ali Wehbe=2

\item 4 citations by Jaime Edilberto Munoz Rivera, collaboration distance with Abdelaziz Soufyane=1, collaboration distance with Ali Wehbe=2

\item 3 citations by Mounir Afilal, collaboration distance with Abdelaziz Soufyane=1, collaboration distance with Ali Wehbe=2

\item 3 citations by Abdelhak Djebabla, collaboration distance with Abdelaziz Soufyane=3, collaboration distance with Ali Wehbe=2

\item 3 citations by Mirelson M. Freitas, collaboration distance with Abdelaziz Soufyane=1, collaboration distance with Ali Wehbe=2

\item 3 citations by Assane Lo, collaboration distance with Abdelaziz Soufyane=3, collaboration distance with Ali Wehbe=3

\item 3 citations by Denis Mercier, collaboration distance with Abdelaziz Soufyane=2, collaboration distance with Ali Wehbe=1

\item 3 citations by Serge Nicaise, collaboration distance with Abdelaziz Soufyane=2, collaboration distance with Ali Wehbe=1 

\item 3 citations by  Yu Ming Qin, collaboration distance with Abdelaziz Soufyane=2, collaboration distance with Ali Wehbe=2

\item 3 citations by Reinhard Racke, collaboration distance with Abdelaziz Soufyane=2, collaboration distance with Ali Wehbe=3

\item 3 citations by Carlos Alberto Raposo, collaboration distance with Abdelaziz Soufyane=1, collaboration distance with Ali Wehbe=2

\item 3 citations by Xinguang Yang, collaboration distance with Abdelaziz Soufyane=2, collaboration distance with Ali Wehbe=2

\item 2 citations by Mohammad Akil, collaboration distance with Abdelaziz Soufyane=1, collaboration distance with Ali Wehbe=1

\item 2 citations by Moncef Aouadi, collaboration distance with Abdelaziz Soufyane=1, collaboration distance with Ali Wehbe=2

\item 2 citations by Maya Bassam, collaboration distance with Abdelaziz Soufyane=2, collaboration distance with Ali Wehbe=1

\item 2 citations by Amirouche Berkani, collaboration distance with Abdelaziz Soufyane=3, collaboration distance with Ali Wehbe=3

\item 2 citations by Margareth da Silva Alves, collaboration distance with Abdelaziz Soufyane=2, collaboration distance with Ali Wehbe=3

\item 2 citations by Takahiro Endo, collaboration distance with Abdelaziz Soufyane=5, collaboration distance with Ali Wehbe=

\end{itemize}

\begin{oss}\label{AutX}

The above list is incomplete. For instance, Wael Youssef cites \cite{SW2003} twice in publications co-authored with Ali Wehbe (collaboration distance with Abdelaziz Soufyane=2, collaboration distance with Ali Wehbe=1),  and Ammar Khemmoudj cites \cite{SW2003} twice (collaboration distance with Abdelaiz Soufyane=3, collaboration distance with Ali Wehbe=3) in publications co-authored with respectively
Naouel Kechiche and Taklit Hamadouche. Some other authors' citations are not taken into account.

\end{oss}

\subsubsection{Citations of \cite{SW2003} by Institutions}

\begin{itemize}

\item 35 citations by King Fahd University of Petroleum and Minerals

\item 17 citations by Federal University of Par\'a

\item 14 citations by University of Sharjah

\item 7 citations by Department of Mathematical Sciences, King Fahd University of Petroleum and Minerals

\item 9 citations by Universit\'e Badji Mokhtar de Annaba

\item 9 citations by Lebanese (Libanaise) University

\item 8 citations by Southwestern University of Finance and Economics, Chengdu

\item 6 citations by Universit\"at Konstanz

\vdots

\end{itemize}

\begin{oss}\label{Inst}

The above list is incomplete since the list is distributed for a same institution through distinct research centers or departments. Therefore the above data have to be explored more precisely.

\end{oss}

\subsection{Some data extracted from Google Scholar and the open access databasis zbMATH}

The last citations of \cite{SW2003} in 2023 and 2022 are for instance from articles published in the peer-reviewed journals:

\begin{itemize}

\item IEEE Transactions on Automatic Control in 2023

\item Applicable Analysis in 2023

\item Mathematics and Mechanics of Solids in 2023

\item Acta Applicandae Mathematicae in 2023

\item Mathematics in 2023

\item Mathematical Methods in the Applied Sciences in 2023

\item Communications on Applied Mathematics and Computation in 2023

\item SeMA Journal in 2023

\item Kragujevac Journal of Mathematics in 2023, to appear in 2024

\item Mediterranean Journal of Mathematics in 2022

\item Zeitschrift f\"ur angewandte Mathematik und Physik in 2022

\item Evolution Equations and Control Theory in 2022

\item Applicable Analysis twice in 2022

\item Journal of Mathematical Analysis and Applications in 2022

\item Acta Applicandae Mathematicae in 2022 (with self-citation by the first author of \cite{SW2003})

\item Bulletin of the Malaysian Mathematical Sciences Society in 2022

\item International Journal of Computer Mathematics in 2022

\item Asymptotic Analysis in 2021

\item Computational and Applied Mathematics in 2021 (with self-citation by the second author of \cite{SW2003})

\item Discrete and Continuous Dynamical Systems - Series S in 2021

\item Osaka Journal of Mathematics in 2021

\item International Journal of Applied Mathematics and Computer Science in 2021

\item Zeitschrift f\"ur angewandte Mathematik und Physik in 2021

\item Advances in Mathematics: Scientific Journal in 2021

\item Journal of Evolution Equations in 2020

\item Archiv der Mathematik in 2020

\item Acta Mechanica in 2020 (with self-citation by the first author of \cite{SW2003})

\item Asymptotic Analysis in 2020 (with self-citation by the second author of \cite{SW2003})

\item Journal of Mathematical Analysis and Applications in 2020

\item Mathematical Modelling of Natural Phenomena in 2019 

\item Chinese Quarterly Journal of Mathematics in 2019

\item Advances in Nonlinear Analysis in 2018

\item Journal of Mathematical Physics in 2018

\item Journal of Integral Equations and Applications 2018

\item Zeitschrift f\"ur Angewandte Mathematik und Mechanik in 2017

\item Afrika Matematika in 2017

\item Journal of Mathematical Physics in 2017

\item IEEE Transactions on Automatic Control in 2017

\item IMA Journal of Mathematical Control and Information in 2016 (with self-citation by the first author of \cite{SW2003})

\item Discrete and Continuous Dynamical Systems in 2016

\item Acta Mathematica Scientia in 2016

\item Electronic Journal of Differential Equations in 2015

\item Applied Mathematics and Optimization in 2012 and 2015

\item Boundary Value Problems in 2015

\item Journal of Spacecraft and Rockets in 2014

\item Journal of Dynamical and Control Systems twice in 2010, in 2014

\item Differential Integral Equations in 2014

\item Partial Stabilization and Control of Distributed Parameter Systems with Elastic Elements in 2014 (Chapter)

\item Quarterly of Applied Mathematics in 2013

\item International Journal of Analysis and Applications in 2013

\item Journal of Differential Equation in 2012

\item IMA Journal of Applied Mathematics in 2012

\item Journal of Mathematical Analysis and Applications in 2008 2009 and 2011

\item Mathematical Control and  Related Fields in 2011

\item Advances in Differential Equations in 2009

\item Applied Mathematics and Computation in 2008

\item Nonlinear Differential Equations and Applications NoDEA in 2008

\item Journal of Mathematical Physics in 2010, 2011, 2013

\vdots

\end{itemize}

Some self-citations by the second author of \cite{SW2003} are from articles published in the peer-reviewed journals:

\begin{itemize}

\item Journal of Mathematical Analysis and Applications in 2007 (by the first author of \cite{SW2003})

\item Nonlinear Analysis: Real World Applications in 2009 (by the first author of \cite{SW2003})

\item Applicable Analysis in 2009 (by the second author of \cite{SW2003})

\item Journal of Mathematical Physics in 2010 (by the second author of \cite{SW2003})

\item Comptes Rendus de l'Acad\'emie des Sciences in 2011 (by the second author of \cite{SW2003})

 \item Differential Equations \& Applications in 2011 (by the second author of \cite{SW2003})

\item Dynamic Systems and Applications in 2012 (by the second author of \cite{SW2003})

\item Journal of Mathematical Analysis and Applications in 2015 (by the second author of \cite{SW2003})

\item Applicable Analysis in 2017 (by the first author of \cite{SW2003})

\item Comptes Rendus de l'Acad\'emie des Sciences in 2019 by (the second author of \cite{SW2003})

\item Mathematics and Mechanics of Solids in 2020 (by the first author of \cite{SW2003})

\item Acta Mechanica in 2020 (by the first author of \cite{SW2003})

\item Computational and Applied Mathematics in 2021 (by the second author of \cite{SW2003})

\end{itemize}

\vskip 0.2mm

\subsection{Some data extracted from the open access databasis zbMATH}

In the databasis zbMATH \cite{SW2003} is cited 10 times by Salim A. Messaoudi., 7 times by Nasser-eddine Tatar, 6 times by  Muhammad Islam Mustafa, 4 times by Abdulaziz Tijani  Apalara, 4 times by Baowei Feng \ldots

\bigskip

\end{document}